\documentclass[sort&compress,final,usenames, dvipsnames]{elsarticle}
\usepackage{setspace}

\usepackage{amsmath}

\usepackage{soul}

\usepackage{lineno,hyperref}
\usepackage{graphicx}
\modulolinenumbers[5]
\usepackage{float}
\usepackage{subcaption}
\usepackage{caption}
\graphicspath{{./figs/}} 

\journal{Computers \& Fluids}

\usepackage{tikz}
\usepackage{tikz-3dplot}
\usetikzlibrary{calc,positioning,matrix,chains,decorations.pathreplacing,arrows}
\usetikzlibrary{matrix,arrows.meta,quotes,shadows,decorations.pathreplacing,positioning,fadings}
\tdplotsetmaincoords{50}{110}

\begin{document}

\begin{frontmatter}
		
	\title{NPLIC: A Machine Learning Approach to Piecewise Linear Interface Construction}
		
	\author[label1,label2]{Mohammadmehdi Ataei}
	\author[label2]{Markus Bussmann\corref{cor1}}
	\author[label2]{Vahid Shaayegan}
	\author[label3]{Franco Costa}
	\author[label4]{Sejin Han}
	\author[label2]{Chul B. Park}

	\cortext[cor1]{Corresponding author}
	\address[label1]{Vector Institute, 661 University Ave Suite 710, Toronto, ON M5G 1M1, Canada}
	\address[label2]{Department of Mechanical Engineering, University of Toronto, 5 King's College Rd, Toronto, ON M5S 3G8, Canada}
	\address[label3]{Autodesk, Inc., 259-261 Colchester Rd., Kilsyth, VIC. 3137, Australia}
	\address[label4]{Autodesk, Inc. 2353 North Triphammer Rd., Ithaca, NY 14850, USA}

	\begin{abstract}
				
		Volume of fluid (VOF) methods are extensively used to track fluid interfaces in numerical simulations, and many VOF algorithms require that the interface be reconstructed geometrically. For this purpose, the Piecewise Linear Interface Construction (PLIC) technique is most frequently used, which for reasons of geometric complexity can be slow and difficult to implement. Here, we propose an alternative neural network based method called NPLIC to perform PLIC calculations. The model is trained on a large synthetic dataset of PLIC solutions for square, cubic, triangular, and tetrahedral meshes. We show that this data-driven approach results in accurate calculations at a fraction of the usual computational cost, and a single neural network system can be used for interface reconstruction of different mesh types.
				
	\end{abstract}

	\begin{keyword}
		Machine Learning \sep Neural Networks \sep PLIC \sep Piecewise Linear Interface Construction \sep Volume Of Fluid \sep VOF \sep Computational Fluid Dynamics 
	\end{keyword}
		
\end{frontmatter}

\section{Introduction}

In the numerical simulation of multiphase flows, the volume of fluid (VOF) method is widely used to track fluid interfaces through a computational domain (e.g.\ \cite{AULISA2003355,GUEYFFIER1999423,AGBAGLAH2011194,meng, WELCH2000662, KLEEFSMAN2005363, RENARDY2001243}). In this method, a scalar field $\alpha$ denotes the volume fraction of one fluid within each cell. In the case of a liquid-gas system, for example, $\alpha = 1$ in liquid cells, $\alpha = 0$ in gas cells, and $0< \alpha  < 1$ in interface cells. 

Reconstruction of the interface geometry from the $\alpha$ field is an important step for calculating volume fluxes advected across cell boundaries \cite{advection}, and related calculations such as finding the distance between two interfaces \cite{ATAEI2021107698}. Since the pioneering work of Youngs \cite{originalPLIC}, Piecewise Linear Interface Construction (PLIC) has been widely employed to geometrically reconstruct interfaces in VOF simulations. 

\begin{figure}[H]
	\centering
	\includegraphics[width=0.7\linewidth]{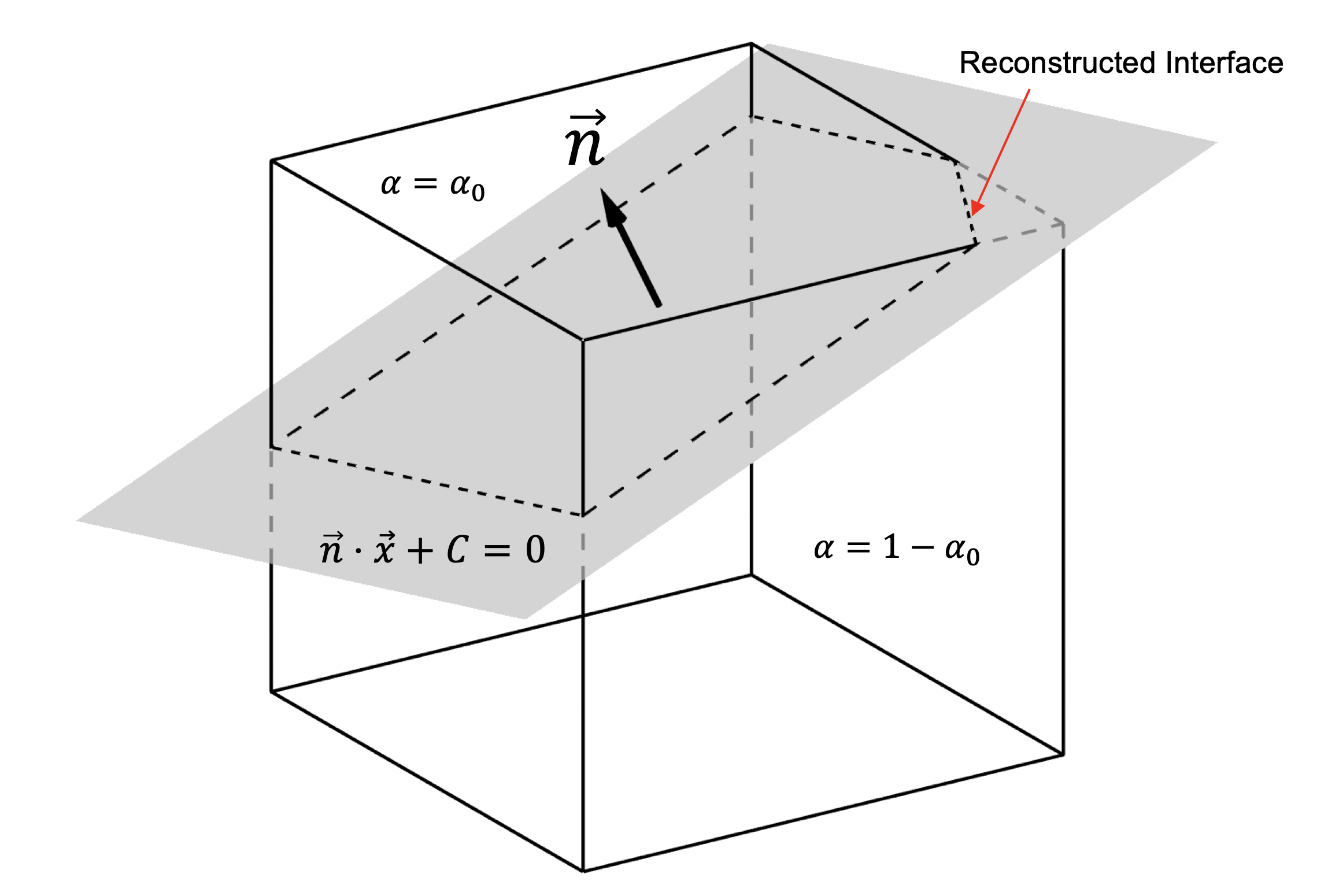}
	\caption{A plane representing the interface splits the cell into two parts, with volume fractions $\alpha_0$ and $1-\alpha_0$.}
	\label{fig:PLIC}
\end{figure}

Given a known interface normal $\vec{n}$ and the volume fraction $\alpha_0$ of an interface cell, PLIC calculates the constant $C$ of the plane $\vec{n}\cdot \vec{x} + C = 0$ ($\vec{x} \in R^{2} \textrm{ in 2D or } R^{3} \textrm{ in 3D}$) that splits the cell into two parts, with volume fractions $\alpha_0$ and $1-\alpha_0$ (see Fig.\ \ref{fig:PLIC}). The reconstructed interface is the polygon resulting from the plane's intersection with the cell.

Especially in 3D, finding $C$ involves complex geometrical operations that can be slow to compute. In some algorithms, $C$ is found iteratively \cite{iterative, iterative2}, until the target volume fractions $\alpha_0$ and $1-\alpha_0$ are achieved within a given tolerance. For simpler geometries such as triangular or rectangular meshes, analytical solutions have been developed that reduce the computational cost \cite{SCARDOVELLI2000228, YANG200641}, although these approaches may also include a slow iterative step used to select one from a set of several governing equations.

Machine learning algorithms are increasingly being applied in computational science \cite{OISHI2017327} and Computational Fluid Dynamics (CFD) simulations in various ways; for example, flow approximation \cite{autodesk, LOPEZPENA2012112, SWISCHUK2019704}, shape optimization for fluid flow processes \cite{Hirschen}, cardiovascular flow modeling \cite{KISSAS2020112623}, shock detection \cite{LIU20191}, computing interface curvature \cite{curvature, PATEL2019104263}, and for turbulence modeling \cite{ling_kurzawski_templeton_2016, SARGHINI200397, ZHOU2019104319}. In this work, we will demonstrate a machine learning approach to find $C$, by using artificial neural networks to relate $C$, $\alpha_0$, and $\vec{n}$ for different cell geometries. It will be shown that neural networks can outperform standard PLIC algorithms, while being nearly as accurate. We will also show that a single neural network system can find the PLIC solution for different mesh types. We limit the results of this paper to square, cubic, triangular, and tetrahedral mesh structures, although the same methodology could be extended to other mesh types.

\section{Methodology}

\begin{figure}[H]
	\centering
	\includegraphics[width=0.7\linewidth]{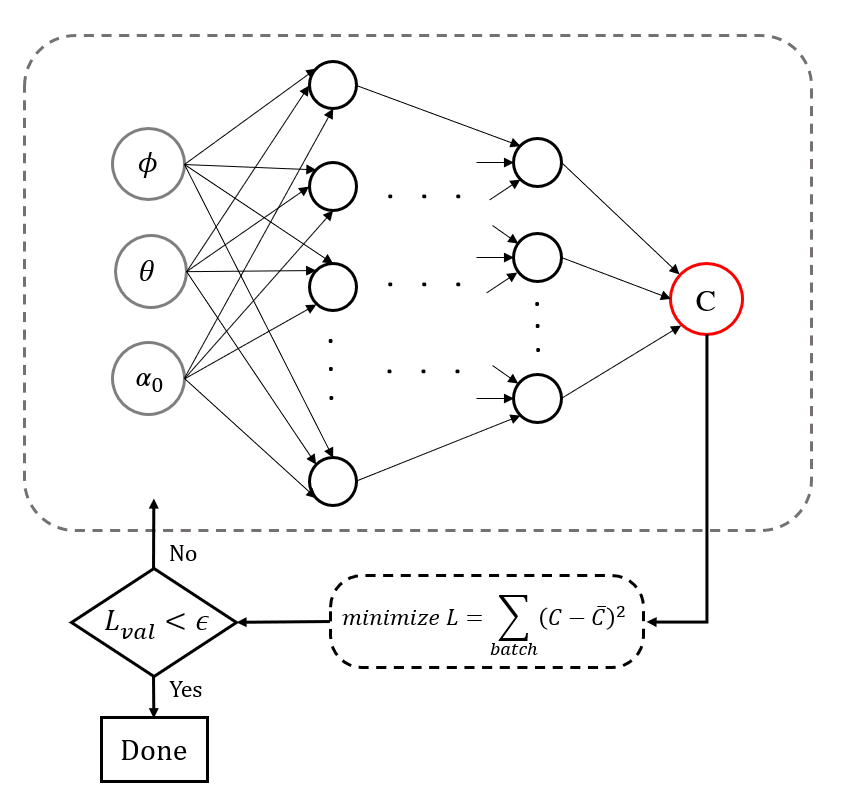}
	\caption{Neural network architecture.}
	\label{fig:NN}
\end{figure}

Artificial neural networks can be thought of as universal approximators capable of extracting nonlinear relationships between different parameters through a kind of machine perception \cite{MLP}. Fig.\ \ref{fig:NN} shows a multilayer perceptron (MLP) ``fully-connected" neural network that consists of an interconnected network of so-called artificial neurons, comprised of an input layer, a series of fully-connected hidden layers, and an output layer. Each neuron is made up of a set of inputs, weights, and a bias. The bias is added to the combined sum of input-weight products, and the result passes though an activation unit, which is usually a sigmoid or ReLU function.

Initially, the network weights and biases are ignorant of the inputs and outputs. The goal is to find a combination of weights and biases that best map the appropriate inputs to the correct PLIC solution. To this reason, the network is fed a synthetic dataset of PLIC solutions, and the weights and biases are updated iteratively by using the gradient of the following loss function:

\begin{equation}
	L(C,\bar{C}) = \sum_{batch} (C - \bar{C})^2
\end{equation}

\noindent which is the squared sum of the difference between each predicted value $C$ and the actual value $\bar{C}$, summed over the training batch. As shown in  Fig.\ \ref{fig:NN}, the training is continued until the validation loss is less than a given tolerance, or a specified maximum number of epochs is reached. 

Initially, we use four separate fully-connected neural networks to find the PLIC constants $C$ for square, cubic, triangular, and tetrahedral meshes. The output layer of each of these networks is a single neuron outputting a constant $C$; the inputs differ depending on the mesh structure. Later on, we will use a single neural network for PLIC calculation of different mesh types.

For square and cubic meshes, the geometries are constant, and so $C$ can be defined as only a function of the interface cell normal $\vec{n}$ and volume fraction $\alpha_0$. As such, for square and cubic meshes the input layers are $(\alpha_0, \theta)$ and $(\alpha_0, \phi, \theta)$ respectively, where $\vec{n} = [cos\ \theta, sin\ \theta]$ in 2D, and $\vec{n} = [sin\ \theta cos\ \phi, sin\ \theta sin\ \phi, cos\ \theta]$ in 3D, in spherical/polar coordinate systems with axes aligned with the cell sides. In VOF models, the interface normal $\vec{n}$ is often obtained from the gradient of the volume fraction field or a smoothed representation of it (e.g., $\vec{n} = \nabla\alpha / |\alpha|$) \cite{normalCalculation}.

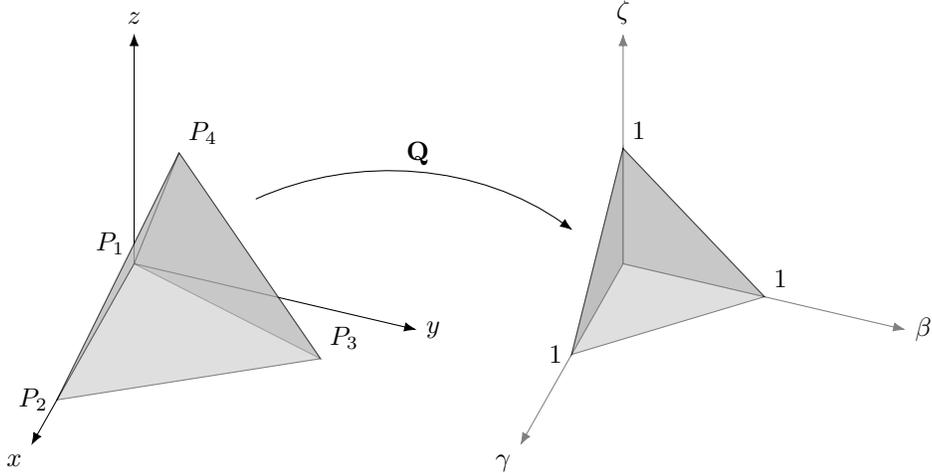
\begin{figure}[H]
\centering	    
	\begin{tikzpicture}[tdplot_main_coords, >=Latex]
		\draw[->] (0,0,0) -- (4,0,0) node[below left] {$x$};
		\draw[->] (0,0,0) -- (0,4,0) node[right] {$y$};
		\draw[->] (0,0,0) -- (0,0,4) node[above] {$z$};
		\begin{scope}[fill opacity=.5, draw opacity=.5, text opacity=1]
			\draw [fill=gray!70] (3,0,0) coordinate [label=left:$P_2$] (D1) -- (1,1,3) coordinate [label=above right:$P_4$] (C1) -- (0,0,0) coordinate [label=above left:$P_1$] (B1) -- cycle;
			\draw [fill=gray!70] (B1) -- (C1) -- (1,3,0) coordinate [label=above right:$P_3$] (A1) -- cycle;
			\draw [fill=gray!50] (D1) -- (C1) -- (A1) -- cycle;
		\end{scope}
		\begin{scope}[xshift=65mm, fill opacity=.5, draw opacity=.5, text opacity=1]
			\draw[->] (0,0,0) -- (4,0,0) node[below left] {$\gamma$};
			\draw[->] (0,0,0) -- (0,4,0) node[right] {$\beta$};
			\draw[->] (0,0,0) -- (0,0,4) node[above] {$\zeta$};
			\draw [fill=gray] (2,0,0) coordinate [label=left:$1$] (D) -- (0,0,2) coordinate [label=above right:$1$] (C) -- (0,0,0) coordinate [label={}] (B) -- cycle;
			\draw [fill=gray!70] (B) -- (C) -- (0,2,0) coordinate [label=above right:$1$] (A) -- cycle;
			\draw [fill=gray!50] (D) -- (C) -- (A) -- cycle;
		\end{scope}
		\draw [->, shorten >=5mm, shorten <=5mm] ($(C1)!.3!(A1)$) [bend left] to node [midway, above, yshift=1pt] {$\mathbf{Q}$} (1,7,3);
	\end{tikzpicture}  
	\caption{Normalizing a tetrahedral cell.}
	\label{fig:normTet}
\end{figure}

As shown in Fig.\ \ref{fig:normTet}, the geometry of an arbitrarily-shaped tetrahedral mesh cell can be expressed by coinciding one of the vertices $P_1$ with the origin, $P_2$ with the axis $x$, and one of the faces with the $x-y$ plane that contains $P_1$, $P_2$, and $P_3$. This tetrahedron can be transformed into a unit tetrahedron in the coordinate system $\vec{\delta}$ using the following transformation $\mathbf{Q}$:

\begin{equation}
\underbrace{
	\begin{pmatrix}
		x \\
		y \\
		z 
	\end{pmatrix}}_{\vec{x}} = 
	\underbrace{\begin{pmatrix}
		x_2 & x_3 & x_4 \\
		0   & y_3 & y_4 \\ 
		0   & z_3 & z_4 
	\end{pmatrix}}_{\mathbf{Q}}
	\underbrace{\begin{pmatrix}
		\gamma \\
		\beta  \\
		\zeta  
	\end{pmatrix}}_{\vec{\delta}}
\end{equation}

The transformation $\mathbf{Q}$ reduces the number of parameters for the PLIC calculation to $\phi$, $\theta$, and $\alpha_0$. While $\mathbf{Q}$ preserves $\alpha_0$, the unit normal undergoes a nonlinear transformation so $\phi$ and $\theta$ must be recalculated under the transformation $\mathbf{Q}$ \cite{kromer2021facebased}.

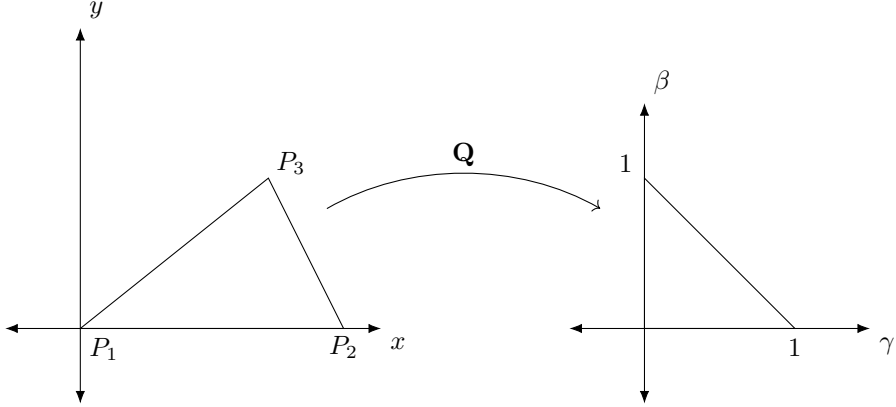
\begin{figure}[H]
\centering
	\begin{tikzpicture}
		\begin{scope}
			    
			\draw [{Latex[]}-{Latex[]}] (-1,0) -- (4,0) node [below right] {$x$};
			\draw [{Latex[]}-{Latex[]}] (0,-1) -- (0,4) node [above right] {$y$};
			\coordinate (A) at (0,0);
			\coordinate (B) at (3.5,0);
			\coordinate (C) at (2.5, 2);
			\draw[draw=black] (A) node [below right] {$P_1$} -- (B) -- (C) -- cycle;
			
			\node at (3.5,-2.5mm){$P_2$};
			\node at (2.8,2.2){$P_3$};
		\end{scope}
		\begin{scope}[xshift=75mm]
			    
			\draw [{Latex[]}-{Latex[]}] (-1,0) -- (3,0) node [below right] {$\gamma$};
			\draw [{Latex[]}-{Latex[]}] (0,-1) -- (0,3) node [above right] {$\beta$};
			\coordinate (A) at (0,0);
			\coordinate (B) at (2,0);
			\coordinate (C) at (0, 2);
            \draw[draw=black] (A) node [below right] {} -- (B) -- (C) -- cycle;
		
			\node at (2,-2.5mm){$1$};
			\node at (-2.5mm,2.2){$1$};
			    
		\end{scope}
		\draw [->, shorten >=5mm, shorten <=5mm] (4,2.5,3) [bend left] to node [midway, above, yshift=1pt] {$\mathbf{Q}$} (8.5,2.5,3) ;
		
	\end{tikzpicture}
	\caption{Normalizing a triangular cell.}
	\label{fig:normTriangular}
\end{figure}

Similarly, as shown in Fig.\ \ref{fig:normTriangular}, an arbitrarily-shaped triangular cell can be transformed into a unit triangle under $\mathbf{Q}$: 

\begin{equation}
	\underbrace{\begin{pmatrix}
		x \\
		y \\
	\end{pmatrix}}_{\vec{x}} = 
	\underbrace{\begin{pmatrix}
		x_2 & x_3 \\
		0   & y_3 
	\end{pmatrix}}_{\mathbf{Q}}
	\underbrace{\begin{pmatrix}
		\gamma \\
		\beta  \\
	\end{pmatrix}}_{\vec{\delta}}
\end{equation}

\noindent In this case, the input layer of the network is $\theta$ and $\alpha_0$.

The fact that the inputs to the neural network for triangular and square meshes in 2D, and the inputs for tetrahedral and cubic meshes in 3D are identical, also allows us to use a single neural network to calculate a PLIC solution for triangular and square meshes, and similarly use another neural network for tetrahedral and cubic mesh types. An input parameter $m$ is introduced to distinguish between the mesh types for the neural networks, which takes the value of 1 for triangular/tetrahedral cell types and zero for square/cubic mesh types. 

In summary, the input parameters for the neural networks depending on the mesh type are given in Table \ref{table:nninputs}:

\begin{table}[H]
	\centering
	\begin{tabular}{c|c}
		Mesh type   & Parameters                 \\ \hline \hline
		Square      & $\theta$, $\alpha_0$         \\
		Triangular  & $\theta$, $\alpha_0$         \\
		Cubic       & $\phi$, $\theta$, $\alpha_0$ \\
		Tetrahedral & $\phi$, $\theta$, $\alpha_0$ \\
		Triangular and square (T-S)   & $\theta$, $\alpha_0$, $m$         \\
		Tetrahedral and cubic (T-C)   &$\phi$, $\theta$, $\alpha_0$, $m$   \\      
	\end{tabular}
	\caption{Inputs of the neural network for different mesh types.}
	\label{table:nninputs}
\end{table}

We generated synthetic PLIC datasets for each mesh type by computing the PLIC solution (i.e., the constant $C$) for a large number of input parameters. The following equations were used to discretely approximate the boundary of a unit sphere to create sets of normal orientations $S_{\vec{n}}$ and volume fractions $S_{\alpha_0}$ (adapted from \cite{maric2020iterative, kromer2021facebased}:

\begin{eqnarray}
    \small
    S_{\vec{n}}(N_{\vec{n}}) = \bigg\{[cos\phi\ sin\theta, sin\phi\ sin\theta, cos \theta]: (\phi, \theta) \in \nonumber \\ \frac{\pi}{2 N_{\vec{n}}}[1, 2,\dots, 2N_{\vec{n}}] \times \frac{\pi}{N_{\vec{n}}}[0, 1,\dots, N_{\vec{n}}] \bigg\},
    \label{eqn:s3d}
\end{eqnarray}
\begin{eqnarray}
    S_{\alpha}(N_{\alpha}) = \{10^{-k}: 5 \leq k \leq 9\} \cup \nonumber \\ \big\{10^{-4} + \frac{m-1}{N_{\alpha} - 1} \big(1 - 2 \cdot 10^{-4}\big): 1 \leq m \leq N_{\alpha} \big\} \cup \nonumber \\ \{1 - 10^{-k}: 5 \leq k \leq 9\}
    \label{eqn:sa}
\end{eqnarray}

\noindent Eq.\ \ref{eqn:sa} ensures that $\alpha_0$ is sufficiently sampled near $0$ and $1$. In 2D, Eq.\ \ref{eqn:s3d} becomes:

\begin{equation}
    S_{\vec{n}}(N_{\vec{n}}) = \bigg\{[cos\theta, sin\theta]: (\theta) \in \frac{\pi}{N_{\vec{n}}}[0, 1,\dots, N_{\vec{n}}] \bigg\}
    \label{eqn:s2d}
\end{equation}

We set $N_{\alpha}=20$ and $N_{\vec{n}}=40$ in 3D, which results in $2N_{\vec{n}}(N_{\vec{n}} + 1)(N_{\alpha} + 10)=98400$ input parameters, and $N_{\alpha}=100$ and $N_{\vec{n}}=100$ in 2D, that gives us $N_{\vec{n}}(N_{\alpha} + 10) = 11000$ input parameters.

Each dataset was randomly split into three parts: $70\%$ for training, $20\%$ for testing, and $10\%$ for validation. Each model was initially trained on the training dataset; during training, the validation dataset was used to prevent over-fitting; and the test dataset was used to evaluate the predictive performance of the final trained model.

The neural network models were developed and trained using the Pytorch \cite{NEURIPS2019_bdbca288} deep learning library. Computations were carried out on an NVIDIA$^{\small{\textregistered}}$ 1080 Ti Graphics Card with 11GB GDDR5X frame buffer, using a PC running Linux Ubuntu 20.04 with Intel$^{\small{\textregistered}}$ Core™ i7-8700K Coffee Lake Processor (6 Cores, up to 4.7 GHz) and 32GB of DDR4 RAM. 

In this work, we use a deep neural network with one hidden layer containing $N$ neurons, and ReLU as the activation function \cite{Hanin_2019}, except for the output layer which is a linear function to allow for negative outputs. Each neural network was trained until the validation loss was less than $\epsilon=5\times 10^{-5}$, up to a maximum of $10^{5}$ epochs with a batch size of 8192, using $N=24$ or $N=48$ in the hidden layer. For training, the Adam optimization algorithm \cite{adam} was used with a learning rate of $10^{-4}$. The training of each network took 10-20 minutes. The results presented in this work are the average of three runs with three different random seeds for the initialization of the networks.

\begin{figure}[H]
	\centering
	\includegraphics[width=1.\textwidth]{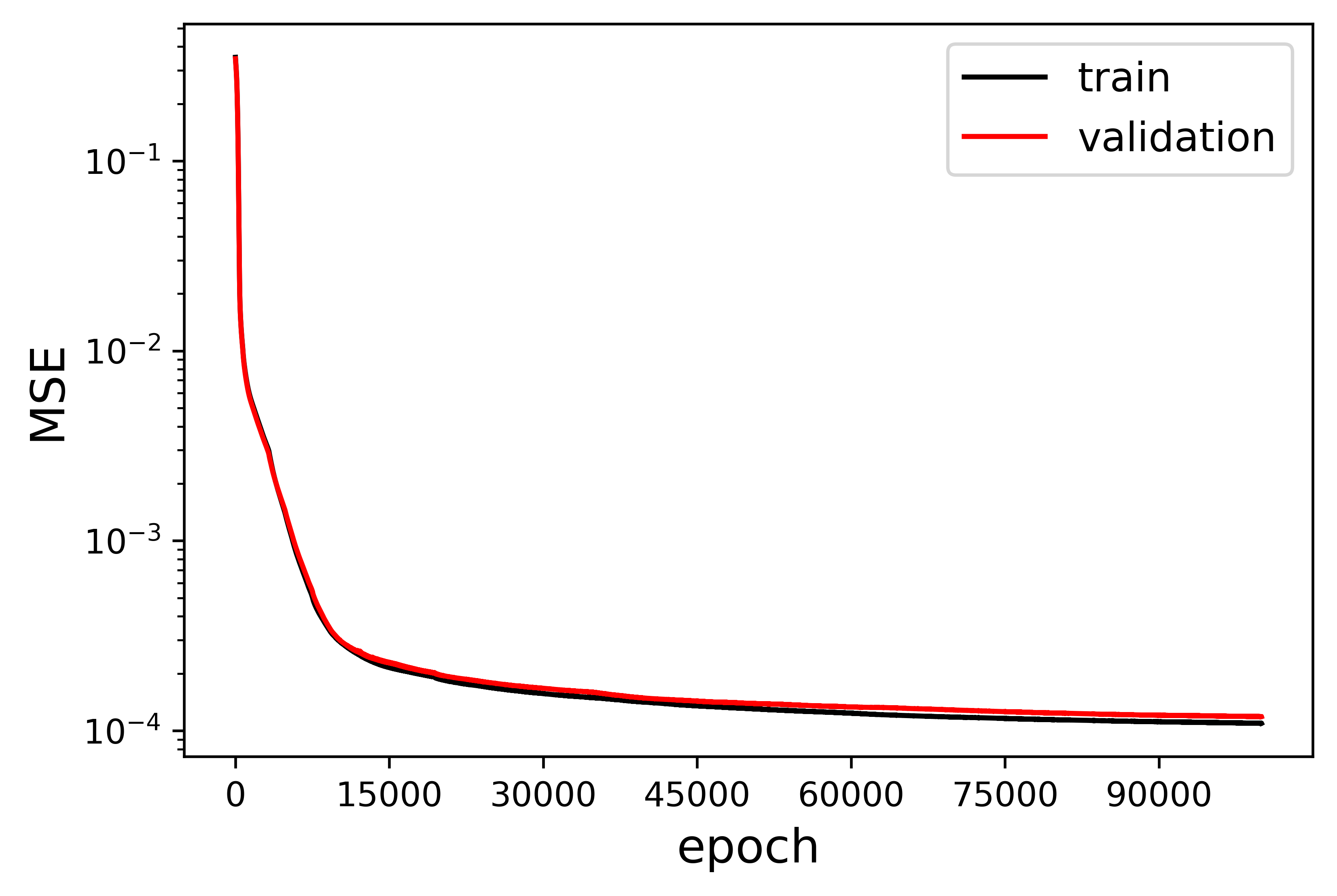}
	\caption{Mean Squared Error (MSE) vs epoch for a neural network with $N=24$ trained on the cubic mesh dataset.}
	\label{training}
\end{figure}

The training performance was evaluated using the Mean Squared Error (MSE). Fig.\ \ref{training} shows MSE of the training dataset after each epoch for a neural network with $N=24$ trained on the cubic mesh dataset.

The neural network models for computing PLIC calculations are hereafter referred to as NPLIC.

\section{Results and Discussion}

\subsection{Predictive Performance}

\begin{figure}[H] 
	\centering
	\begin{minipage}[b]{0.5\linewidth}
		\centering
		\includegraphics[width=\linewidth]{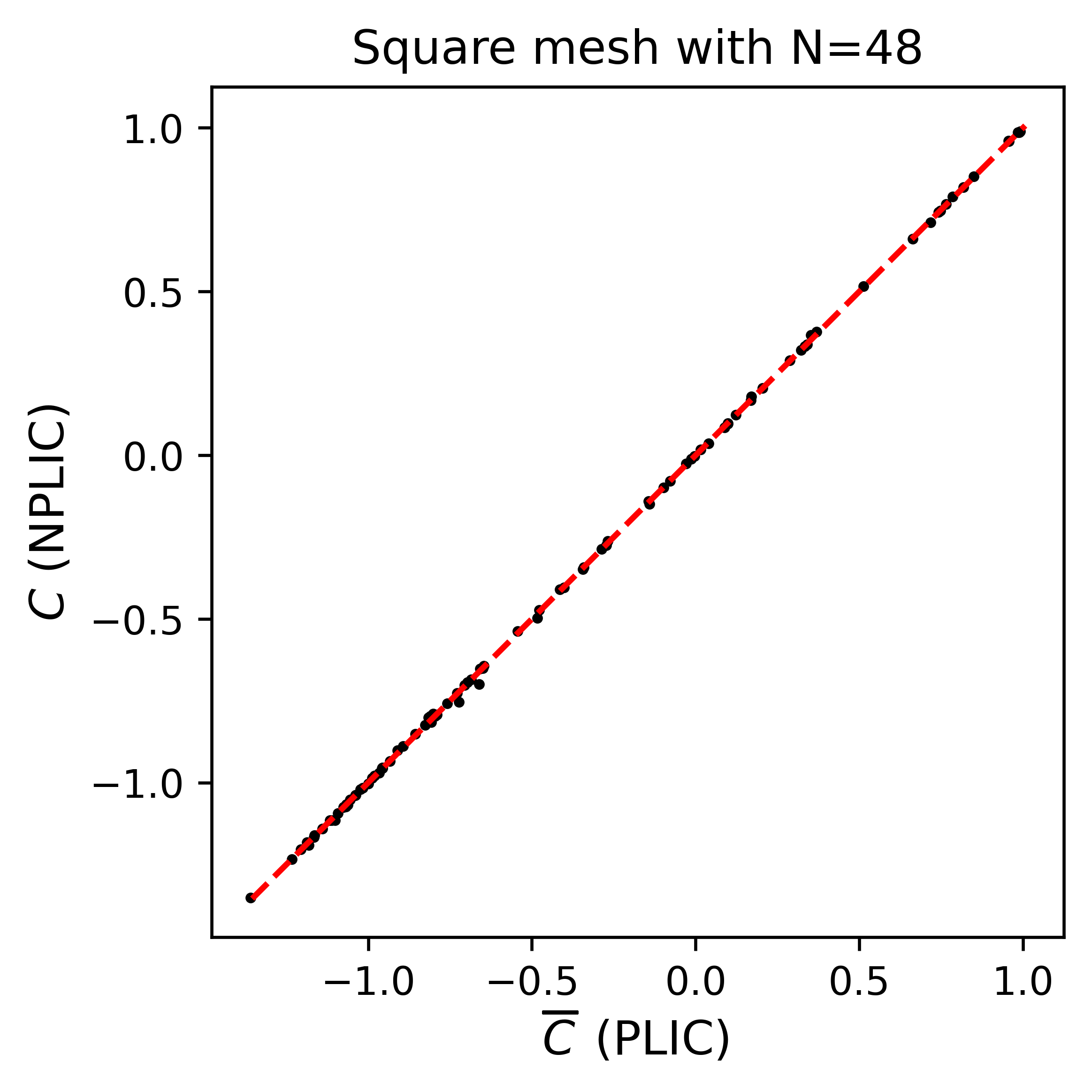}
						
	\end{minipage}
	\begin{minipage}[b]{0.5\linewidth}
		\centering
		\includegraphics[width=\linewidth]{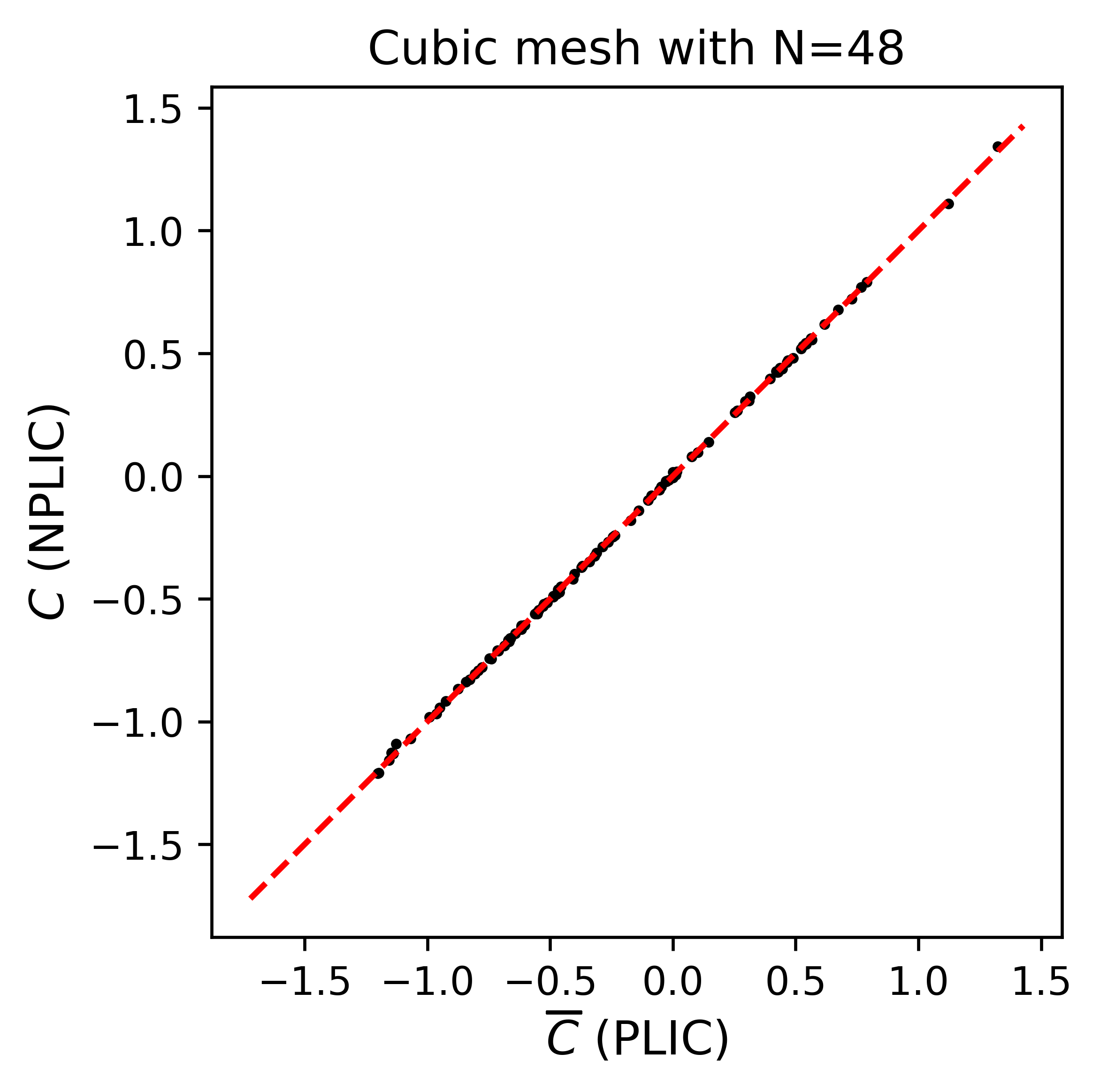} 

	\end{minipage} 
	
	\begin{minipage}[b]{0.5\linewidth}
		\centering
		\includegraphics[width=\linewidth]{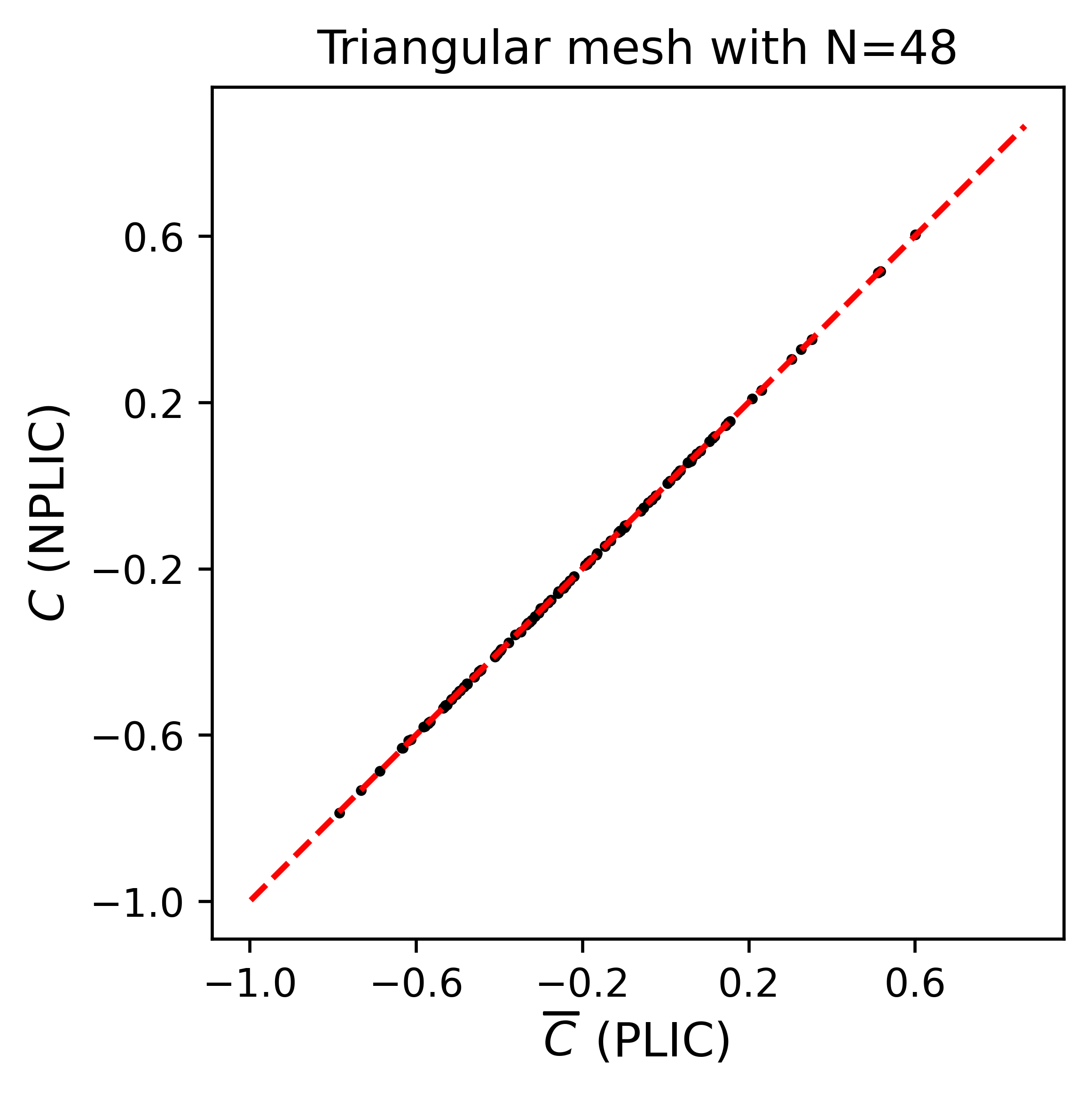}
						
	\end{minipage}
	\begin{minipage}[b]{0.5\linewidth}
		\centering
		\includegraphics[width=\linewidth]{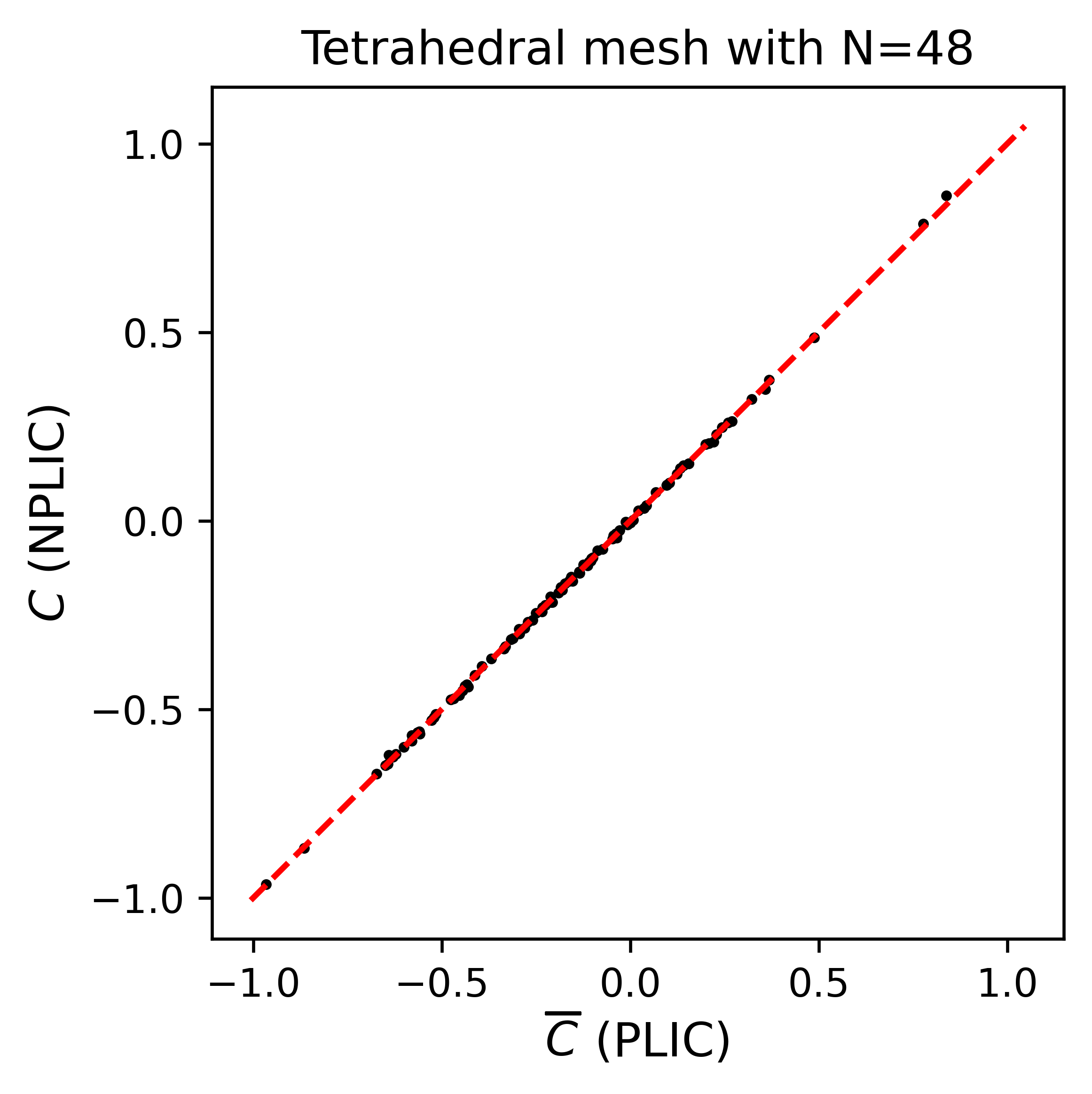} 

	\end{minipage} 
	\begin{minipage}[b]{0.5\linewidth}
		\centering
		\includegraphics[width=\linewidth]{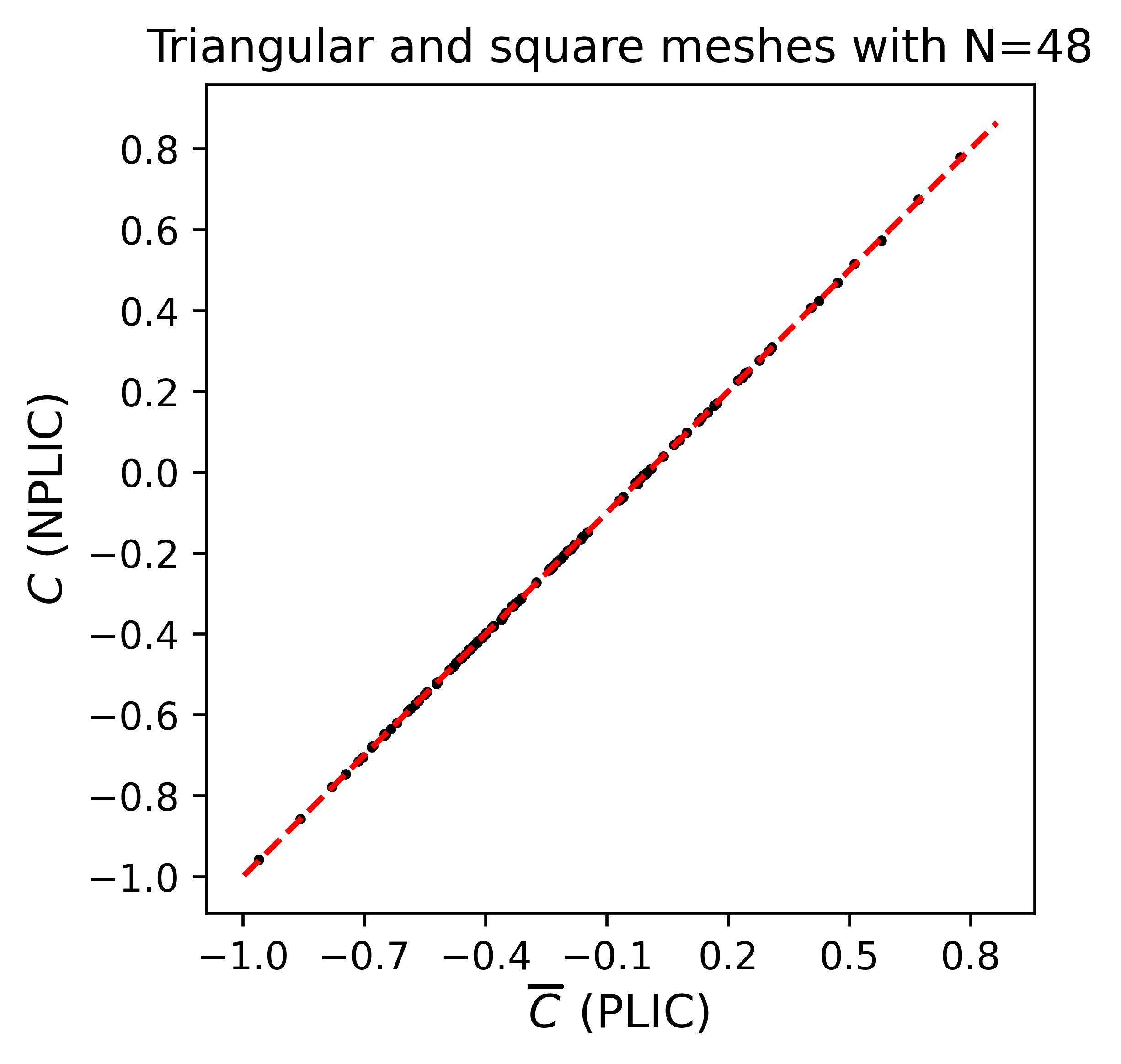}
						
	\end{minipage}
	\begin{minipage}[b]{0.5\linewidth}
		\centering
		\includegraphics[width=\linewidth]{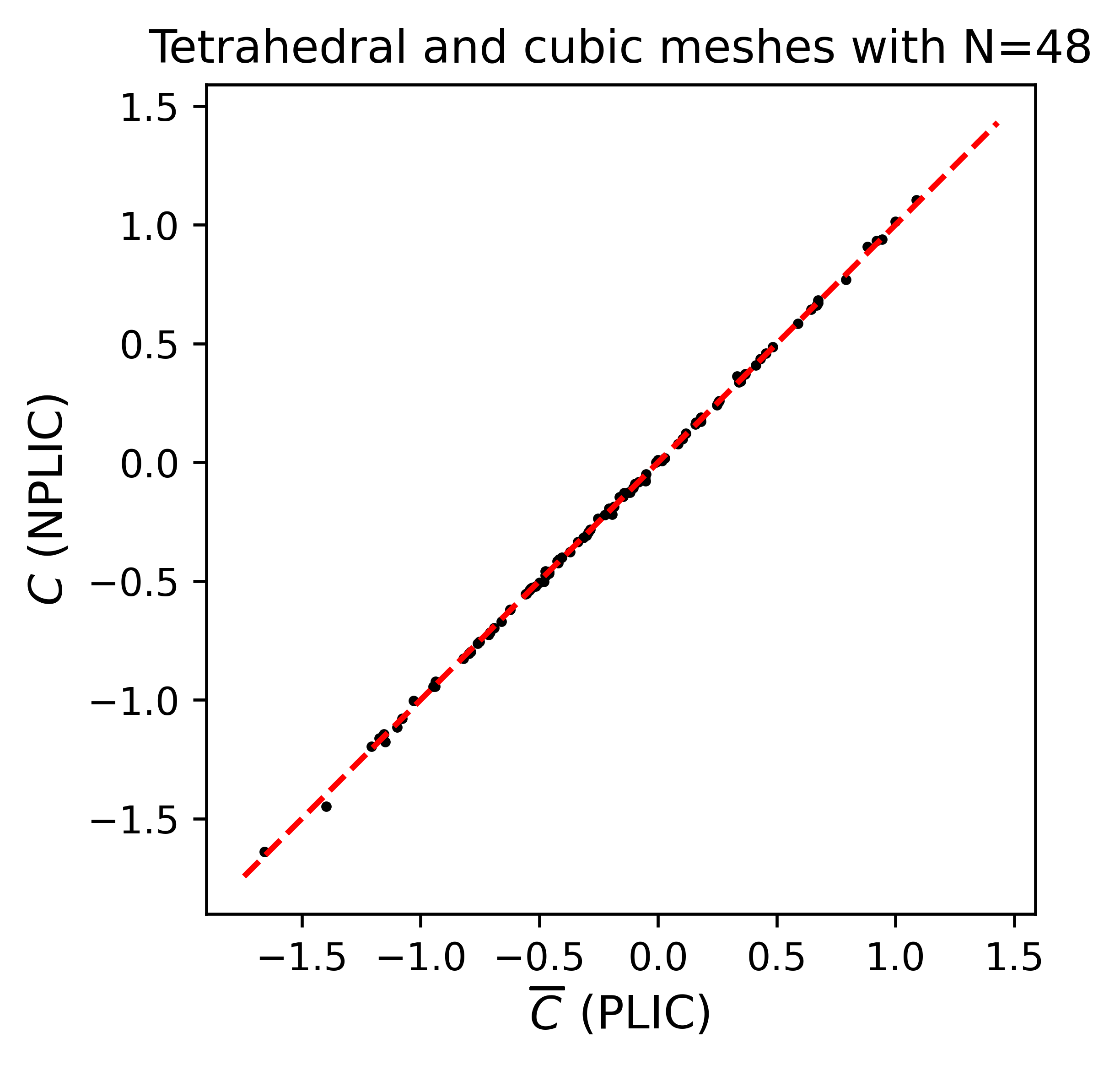} 

	\end{minipage} 
					
	\caption{Plots of the NPLIC predictions vs the test data for each mesh. On each plot, only a hundred data points (selected randomly) are plotted to avoid clutter.}
	\label{prediction}
\end{figure}

\begin{figure}[H]
    \centering
    \includegraphics[width=\linewidth]{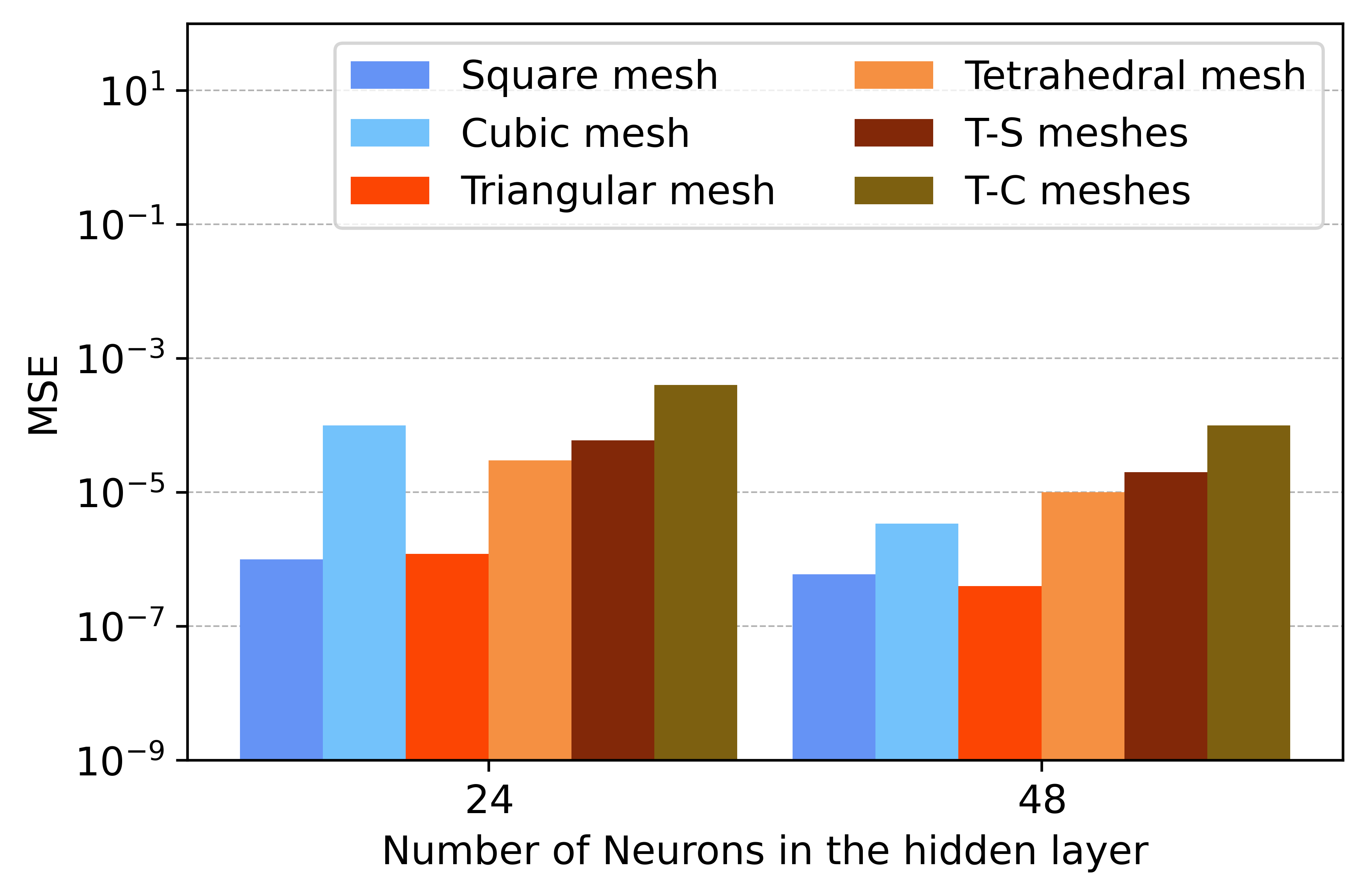}
    \caption{MSE error of NPLIC over the test dataset.}
    \label{fig:MSE}
\end{figure}

\begin{figure}[H]
    \centering
    \includegraphics[width=\linewidth]{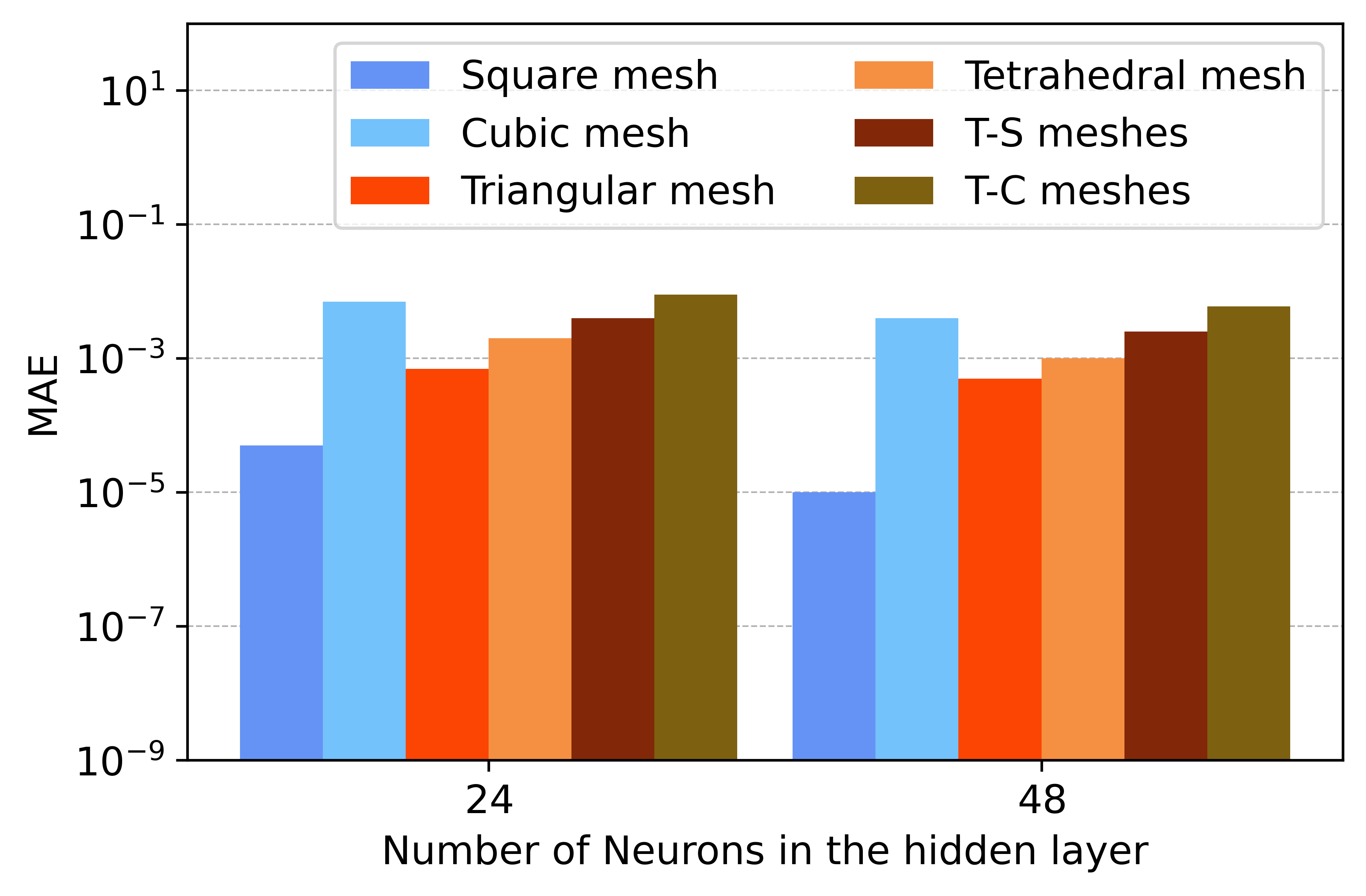}
    \caption{MAE error of NPLIC over the test dataset.}
    \label{fig:MAE}
\end{figure}

The test datasets were used to evaluate the predictive capability of each trained model. In Fig.\ \ref{prediction}, the value of $C$ computed by NPLIC is plotted against $\bar{C}$ from the test dataset. For a perfect fit, all the points would lie on the diagonal. It can be observed that the values of $C$ predicted by NPLIC are in very good agreement with the test results for all the mesh types. The final MSE and MAE (Mean Absolute Error) for each trained model (over the test dataset) are shown in Figs.\ \ref{fig:MSE} and \ref{fig:MAE} for the neural networks that were trained up to $10^5$ epochs. The models with $N=48$ have errors close to $0.1\%$, which implies that NPLIC can be used in place of PLIC algorithms with little to no effect on accuracy. In general, the errors are comparable to other PLIC approximation techniques (e.g.\ \cite{KAWANO2016130}). With increasing mesh complexity, predictably the accuracy of NPLIC decreases, and by increasing the number of neurons the accuracy increases.

\subsection{Speedup}

\begin{figure}[H]
    \centering
    \includegraphics[width=\linewidth]{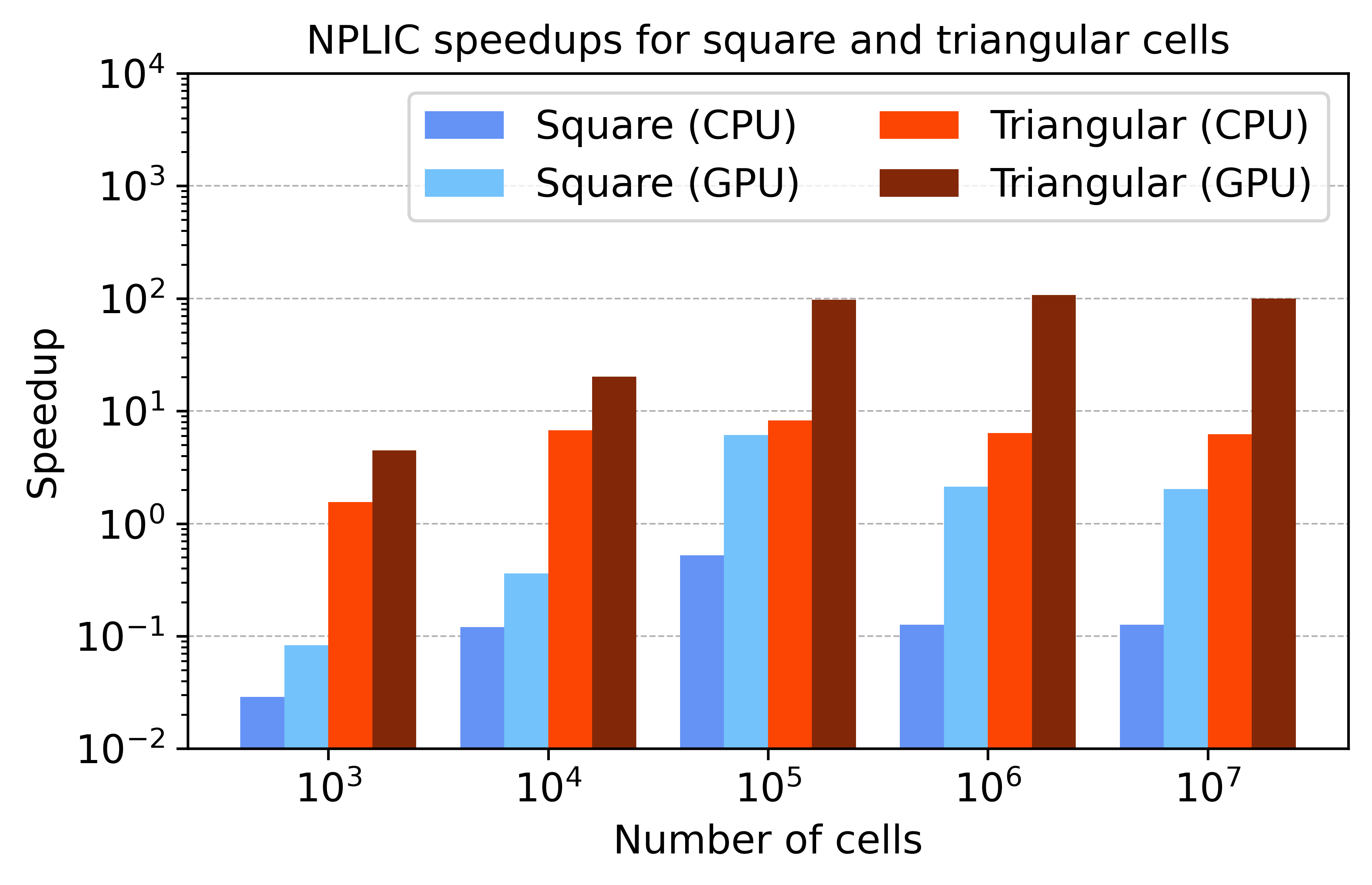}
    \caption{Performance speedups of NPLIC for each 2D mesh type, in comparison to the methods of Scardovelli et al.\ \cite{SCARDOVELLI2000228} (square mesh, analytical) and L\'{o}pez et al.\ \cite{voftools} (triangular mesh, analytical).}
    \label{fig:2Dspeedup}
\end{figure}

\begin{figure}[H]
    \centering
    \includegraphics[width=\linewidth]{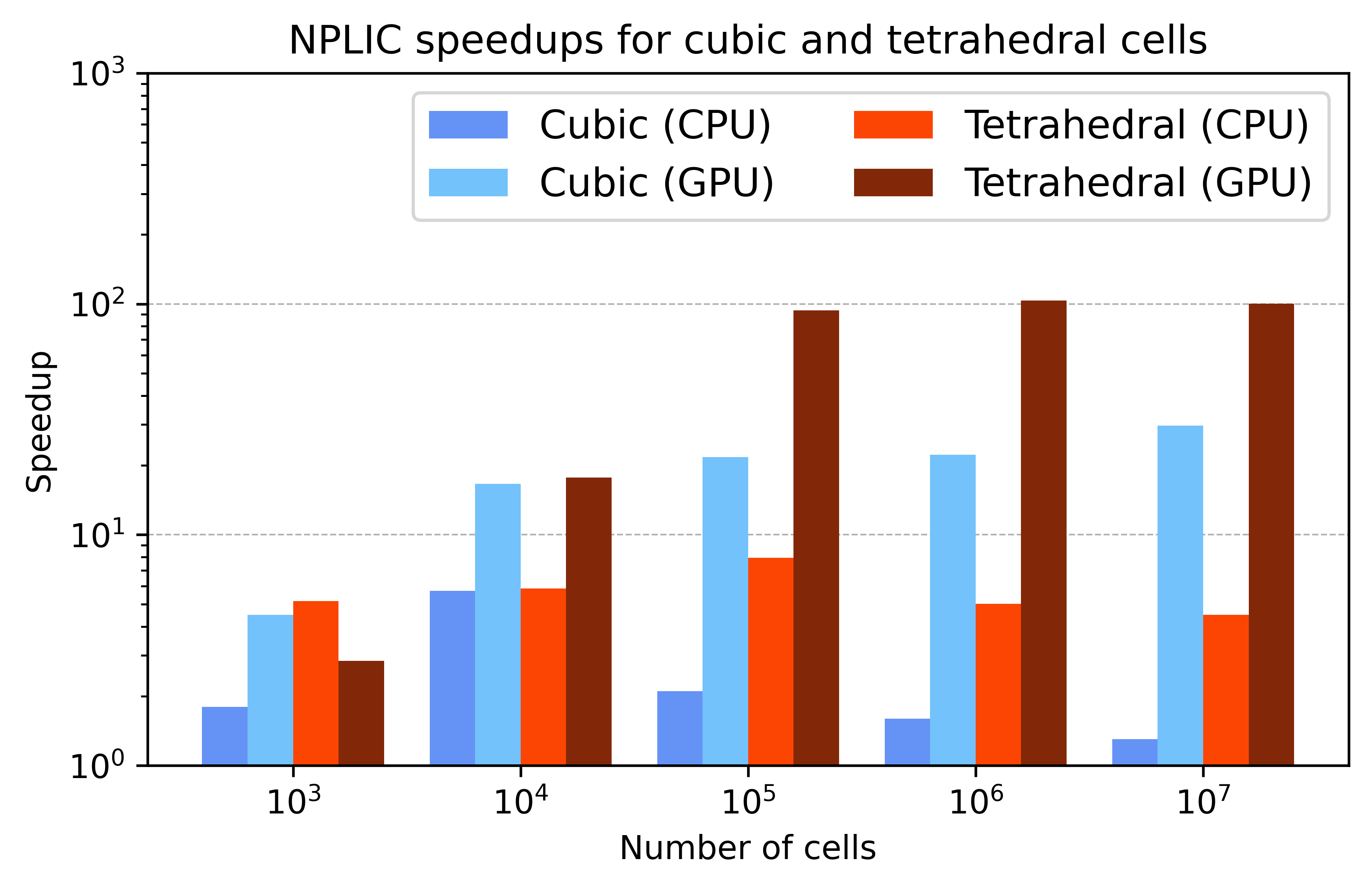}
    \caption{Performance speedups of NPLIC for each 3D mesh type, in comparison to the methods of Scardovelli et al.\ \cite{SCARDOVELLI2000228} (cubic mesh, analytical) and L\'{o}pez et al.\ \cite{voftools} (tetrahedral mesh, analytical).}
    \label{fig:3Dspeedup}
\end{figure}

Figs.\ \ref{fig:2Dspeedup} and \ref{fig:3Dspeedup} show a comparison of NPLIC speedups versus a number of popular PLIC models. To calculate the speedups, the PLIC algorithms were implemented in C++, and compiled using optimization flag -O3 using the GCC compiler v 9.3. Since a neural network can leverage performing calculations on batches of inputs (rather than each cell one by one), we compared the speedups by performing the calculations on different numbers of meshes from $10^3$ to $10^7$. The PLIC models were executed on one CPU core. The NPLIC models were executed on both one CPU core and on the GPU. Since the implementation of PLIC models on the GPU is not trivial, we used the performance of PLIC on the CPU as the basis of our comparisons. 

The wall-time of the PLIC algorithms were measured by the \texttt{std::chrono} library, and the wall-time of NPLIC on the CPU and GPU (ignoring GPU warm-up time) were measured by Pytorch's Profiler. For the profiling tasks, we used the more accurate neural networks with $N=48$.

For square and cubic meshes, we compared NPLIC with the analytical PLIC model of Scardovelli and Zaleski \cite{SCARDOVELLI2000228}, which is faster than iterative methods for rectangular meshes. For triangular and tetrahedral meshes, we compared NPLIC with the efficient analytical model of L\'{o}pez et al.\ \cite{voftools}.

In comparison to the PLIC models, NPLIC is up to 100 times faster when executed on a GPU, and it can be more than 8 times faster on one CPU core. It can be seen that the speedup gains are more significant on the more complex 3D meshes and triangular meshes. However, on the square mesh, the analytical implementation of PLIC is sufficiently simple that it cannot be outperformed by NPLIC. 

To explain the NPLIC speedup gains, we compared the floating point operations (FLOPs) of NPLIC versus L\'{o}pez et al.\ \cite{voftools} for performing the calculations on 1 million tetrahedral cells. For the C++ implementation of the PLIC model, the number of FLOPs was estimated by the Intel's Software Development Emulator (Intel SDE) tool. The Thop library was used to calculate the number of FLOPs for the NPLIC model. We found that while the PLIC model requires fewer floating operations than the NPLIC (1.98 GFLOPs vs 9.6 GFLOPs), on the CPU the NPLIC algorithm operates at higher GFLOPs per second than PLIC (72 GFLOPs per second vs 2.42 GFLOPs per second), and the performance increases further to 1.21 TFLOPs per second on the GPU. This is because NPLIC makes better use of hardware (on both the CPU and GPU), because the essence of neural network operations consists only of simple matrix multiplications, and because deep learning libraries such as Pytorch are heavily optimized to leverage the available hardware.

\subsection{Reconstructing a VOF Scalar Field}
\begin{figure}[H]
	\centering
	\includegraphics[width=0.6\linewidth]{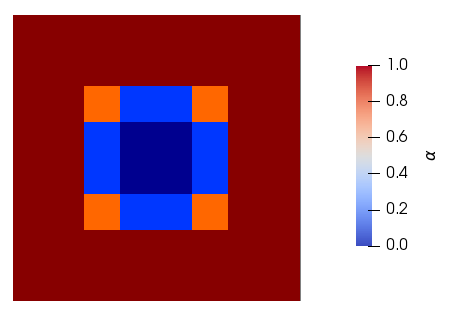}
	\caption{The $\alpha$ scalar field representing a circular bubble on an 8 by 8 grid.}
	\label{voffield}
\end{figure}

\begin{figure}[H]
	\centering
	\includegraphics[width=\linewidth]{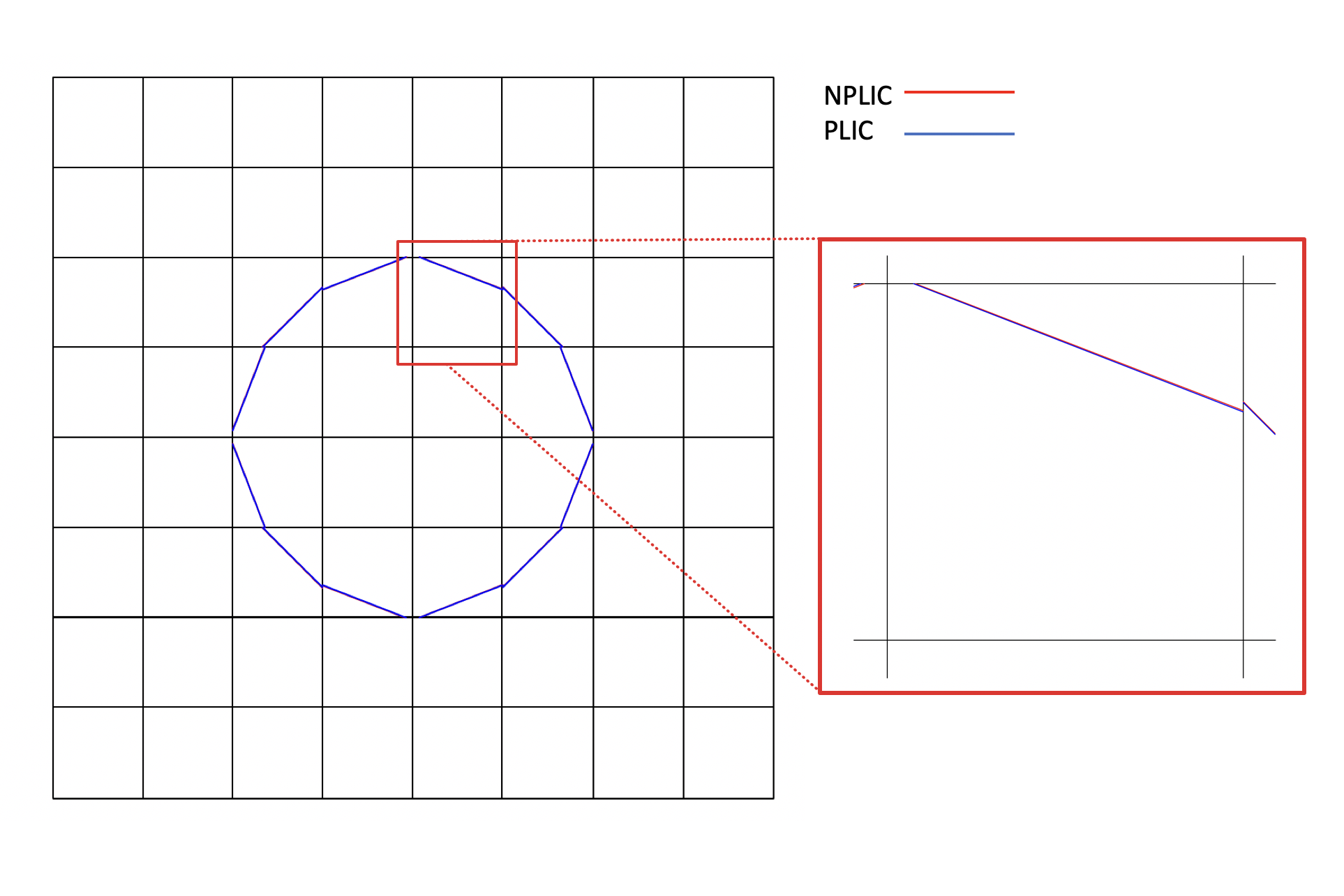}
	\caption{NPLIC and PLIC reconstructions in the bubble interface cells.}
	\label{reconstruction}
\end{figure}

Fig.\ \ref{voffield} shows a 2D VOF scalar field for a circular bubble in a liquid. The gas and liquid phases correspond to $\alpha =0$ and $\alpha = 1$, and the interface cells $0 < \alpha < 1$. A coarse $8\times8$ grid is chosen so that we can easily illustrate the reconstruction results. In Fig.\ \ref{reconstruction}, the interface cells are reconstructed by the Scardovelli analytical PLIC and NPLIC models. The zoomed-in plot shows that the results of NPLIC and PLIC overlap as expected.

\subsection{Implementation in a CFD Solver}

Finally, Basilisk (https://basilisk.fr) \cite{basilisk} is an open-source CFD solver for adaptive Cartesian meshes. Basilisk uses PLIC VOF to solve the advection equation for the volume fraction field \cite{basilisk,advection}, by reconstructing interfaces using PLIC, and calculating volume fluxes geometrically. 

\begin{figure}[H]
	\includegraphics[width=\linewidth]{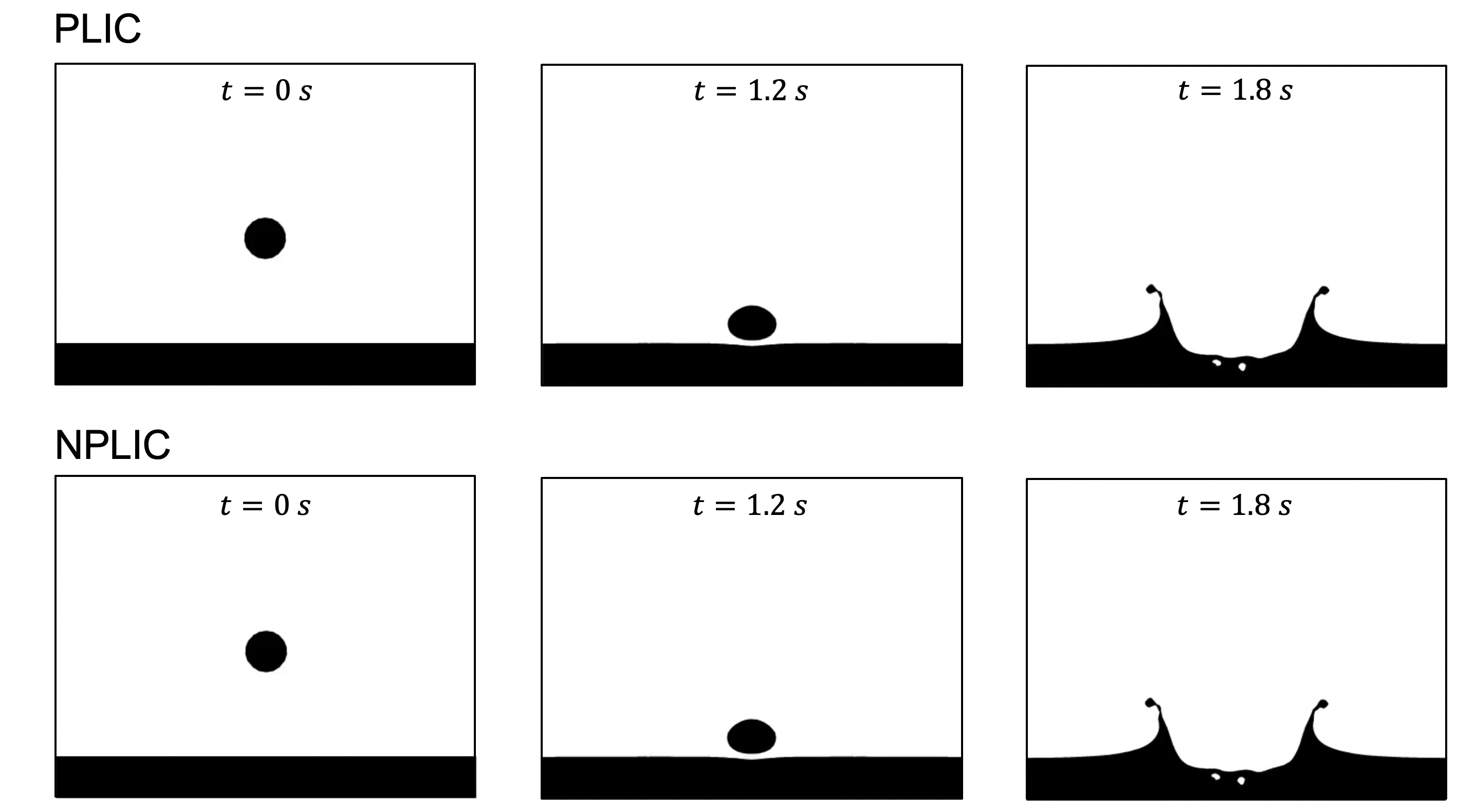}
	\caption{NPLIC performs as well as PLIC when implemented in a multiphase flow solver.}
	\label{splash}
\end{figure}

NPLIC can be easily integrated into an existing codebase. The trained neural networks can be called as subroutines from other programs (e.g.\ using Pytorch's TorchScript for C++ implementations). We replaced the PLIC algorithm in Basilisk with NPLIC, and ran a 2D simulation of a droplet impacting a pool of liquid. As shown in Fig.\ \ref{splash}, the NPLIC results are indistinguishable from the PLIC ones at all timesteps. In this case, NPLIC performed the interface reconstructions about five times faster.

\section{Conclusions}

We have presented a machine learning approach to perform Piecewise Linear Interface Construction (PLIC) on square, cubic, and arbitrarily-shaped triangular and tetrahedral meshes. Each mesh type was normalized to reduce the number of inputs to the neural network. Fully-connected deep neural networks were trained on synthetic datasets for each mesh type. We have shown that the neural networks are capable of performing up to 100 times faster than available PLIC algorithms with minimal loss of accuracy.

\section*{Acknowledgements}

We thank the Natural Sciences and Engineering Research Council of Canada (NSERC) and Autodesk Inc.\ for their financial support.

\bibliography{main}

\begin{thebibliography}{10}
\expandafter\ifx\csname url\endcsname\relax
  \def\url#1{\texttt{#1}}\fi
\expandafter\ifx\csname urlprefix\endcsname\relax\def\urlprefix{URL }\fi
\expandafter\ifx\csname href\endcsname\relax
  \def\href#1#2{#2} \def\path#1{#1}\fi

\bibitem{AULISA2003355}
E.~Aulisa, S.~Manservisi, R.~Scardovelli, S.~Zaleski, A geometrical
  area-preserving volume-of-fluid advection method, Journal of Computational
  Physics 192~(1) (2003) 355--364.

\bibitem{GUEYFFIER1999423}
D.~Gueyffier, J.~Li, A.~Nadim, R.~Scardovelli, S.~Zaleski, Volume-of-fluid
  interface tracking with smoothed surface stress methods for three-dimensional
  flows, Journal of Computational Physics 152~(2) (1999) 423--456.

\bibitem{AGBAGLAH2011194}
G.~Agbaglah, S.~Delaux, D.~Fuster, J.~Hoepffner, C.~Josserand, S.~Popinet,
  P.~Ray, R.~Scardovelli, S.~Zaleski, Parallel simulation of multiphase flows
  using octree adaptivity and the volume-of-fluid method, Comptes Rendus
  Mécanique 339~(2) (2011) 194--207.

\bibitem{meng}
M.~Huang, L.~Wu, B.~Chen, A piecewise linear interface-capturing
  volume-of-fluid method based on unstructured grids, Numerical Heat Transfer,
  Part B: Fundamentals 61~(5) (2012) 412--437.

\bibitem{WELCH2000662}
S.~W. Welch, J.~Wilson, A volume of fluid based method for fluid flows with
  phase change, Journal of Computational Physics 160~(2) (2000) 662--682.

\bibitem{KLEEFSMAN2005363}
K.~Kleefsman, G.~Fekken, A.~Veldman, B.~Iwanowski, B.~Buchner, A
  volume-of-fluid based simulation method for wave impact problems, Journal of
  Computational Physics 206~(1) (2005) 363--393.

\bibitem{RENARDY2001243}
M.~Renardy, Y.~Renardy, J.~Li, Numerical simulation of moving contact line
  problems using a volume-of-fluid method, Journal of Computational Physics
  171~(1) (2001) 243--263.

\bibitem{advection}
R.~Scardovelli, S.~Zaleski, Direct numerical simulation of free-surface and
  interfacial flow, Annual Review of Fluid Mechanics 31~(1) (1999) 567--603.

\bibitem{ATAEI2021107698}
M.~Ataei, V.~Shaayegan, F.~Costa, S.~Han, C.~B. Park, M.~Bussmann, Lbfoam: An
  open-source software package for the simulation of foaming using the lattice
  boltzmann method, Computer Physics Communications 259 (2021) 107698.
\newblock \href {https://doi.org/https://doi.org/10.1016/j.cpc.2020.107698}
  {\path{doi:https://doi.org/10.1016/j.cpc.2020.107698}}.

\bibitem{originalPLIC}
D.~L. Youngs, Time-dependent multi-material flow with large fluid distortion,
  Numerical Methods in Fluid Dynamics (1982) 273--486.

\bibitem{iterative}
M.~Skarysz, A.~Garmory, M.~Dianat, An iterative interface reconstruction method
  for {PLIC} in general convex grids as part of a coupled level set volume of
  fluid solver, Journal of Computational Physics 368 (2018) 254--276.

\bibitem{iterative2}
W.~J. Rider, D.~B. Kothe, Reconstructing volume tracking, Journal of
  Computational Physics 141~(2) (1998) 112--152.

\bibitem{SCARDOVELLI2000228}
R.~Scardovelli, S.~Zaleski, Analytical relations connecting linear interfaces
  and volume fractions in rectangular grids, Journal of Computational Physics
  164~(1) (2000) 228--237.

\bibitem{YANG200641}
X.~Yang, A.~J. James, Analytic relations for reconstructing piecewise linear
  interfaces in triangular and tetrahedral grids, Journal of Computational
  Physics 214~(1) (2006) 41--54.

\bibitem{OISHI2017327}
A.~Oishi, G.~Yagawa, Computational mechanics enhanced by deep learning,
  Computer Methods in Applied Mechanics and Engineering 327 (2017) 327--351.

\bibitem{autodesk}
X.~Guo, W.~Li, F.~Iorio, Convolutional neural networks for steady flow
  approximation, Proceedings of the 22nd ACM SIGKDD International Conference on
  Knowledge Discovery and Data Mining (2016) 481--490.

\bibitem{LOPEZPENA2012112}
F.~L. Peña, V.~D. Casás, A.~Gosset, R.~Duro, A surrogate method based on the
  enhancement of low fidelity computational fluid dynamics approximations by
  artificial neural networks, Computers \& Fluids 58 (2012) 112 -- 119.

\bibitem{SWISCHUK2019704}
R.~Swischuk, L.~Mainini, B.~Peherstorfer, K.~Willcox, Projection-based model
  reduction: {Formulations} for physics-based machine learning, Computers \&
  Fluids 179 (2019) 704 -- 717.

\bibitem{Hirschen}
K.~Hirschen, M.~Schäfer, Bayesian regularization neural networks for
  optimizing fluid flow processes, Computer Methods in Applied Mechanics and
  Engineering 195~(7) (2006) 481--500.

\bibitem{KISSAS2020112623}
G.~Kissas, Y.~Yang, E.~Hwuang, W.~R. Witschey, J.~A. Detre, P.~Perdikaris,
  Machine learning in cardiovascular flows modeling: {Predicting} arterial
  blood pressure from non-invasive {4D} flow {MRI} data using physics-informed
  neural networks, Computer Methods in Applied Mechanics and Engineering 358
  (2020) 112623.

\bibitem{LIU20191}
Y.~Liu, Y.~Lu, Y.~Wang, D.~Sun, L.~Deng, F.~Wang, Y.~Lei, A {CNN}-based shock
  detection method in flow visualization, Computers \& Fluids 184 (2019) 1--9.

\bibitem{curvature}
Y.~Qi, J.~Lu, R.~Scardovelli, S.~Zaleski, G.~Tryggvason, Computing curvature
  for volume of fluid methods using machine learning, Journal of Computational
  Physics 377 (2019) 155--161.

\bibitem{PATEL2019104263}
H.~Patel, A.~Panda, J.~Kuipers, E.~Peters, Computing interface curvature from
  volume fractions: {A} machine learning approach, Computers \& Fluids 193
  (2019) 104263.

\bibitem{ling_kurzawski_templeton_2016}
J.~Ling, A.~Kurzawski, J.~Templeton, Reynolds averaged turbulence modelling
  using deep neural networks with embedded invariance, Journal of Fluid
  Mechanics 807 (2016) 155–166.

\bibitem{SARGHINI200397}
F.~Sarghini, G.~de~Felice, S.~Santini, Neural networks based subgrid scale
  modeling in large eddy simulations, Computers \& Fluids 32~(1) (2003)
  97--108.

\bibitem{ZHOU2019104319}
Z.~Zhou, G.~He, S.~Wang, G.~Jin, Subgrid-scale model for large-eddy simulation
  of isotropic turbulent flows using an artificial neural network, Computers \&
  Fluids 195 (2019) 104319.

\bibitem{MLP}
M.-C. Popescu, V.~E. Balas, L.~Perescu-Popescu, N.~Mastorakis, Multilayer
  perceptron and neural networks, WSEAS Transactions on Circuits and Systems
  8~(7) (2009) 579--588.

\bibitem{normalCalculation}
D.~B. Kothe, R.~C. Mjolsness, {RIPPLE} - a new model for incompressible flows
  with free surfaces, AIAA Journal 30~(11) (1992) 2694--2700.

\bibitem{kromer2021facebased}
J.~Kromer, D.~Bothe, Face-based {Volume-of-Fluid} interface positioning in
  arbitrary polyhedra (2021).

\bibitem{maric2020iterative}
T.~Maric, Iterative {Volume-of-Fluid} interface positioning in general
  polyhedrons with {Consecutive Cubic Spline interpolation} (2020).

\bibitem{NEURIPS2019_bdbca288}
A.~Paszke, S.~Gross, F.~Massa, A.~Lerer, J.~Bradbury, G.~Chanan, T.~Killeen,
  Z.~Lin, N.~Gimelshein, L.~Antiga, A.~Desmaison, A.~Kopf, E.~Yang, Z.~DeVito,
  M.~Raison, A.~Tejani, S.~Chilamkurthy, B.~Steiner, L.~Fang, J.~Bai,
  S.~Chintala, Pytorch: An imperative style, high-performance deep learning
  library, in: Advances in Neural Information Processing Systems, Vol.~32,
  Curran Associates, Inc., 2019, pp. 8026--8037.

\bibitem{Hanin_2019}
B.~Hanin, Universal function approximation by deep neural nets with bounded
  width and {ReLU} activations, Mathematics 7~(10) (2019) 992.

\bibitem{adam}
D.~P. {Kingma}, J.~{Ba}, {Adam: A Method for Stochastic Optimization}, arXiv
  e-prints (2014) arXiv:1412.6980.

\bibitem{KAWANO2016130}
A.~Kawano, A simple volume-of-fluid reconstruction method for three-dimensional
  two-phase flows, Computers \& Fluids 134-135 (2016) 130--145.

\bibitem{voftools}
J.~L\'{o}pez, J.~Hern\'{a}ndez, P.~G\'{o}mez, F.~Faura, {VOFTools - A software
  package of calculation tools for volume of fluid methods using general convex
  grids}, Computer Physics Communications 223 (2018) 45--54.

\bibitem{basilisk}
S.~Popinet, An accurate adaptive solver for surface-tension-driven interfacial
  flows, Journal of Computational Physics 228~(16) (2009) 5838--5866.

\end{thebibliography}

\end{document}